\documentclass[a4paper,12pt]{article}
\usepackage{graphicx}
\usepackage[english]{babel}
\usepackage[latin1]{inputenc}
\usepackage{amsmath}
\usepackage{amssymb}
\usepackage{amsthm}
\usepackage{amsfonts}
\usepackage{graphicx}
\theoremstyle{plain}
\newtheorem{theorem}{Theorem}[section]
\newtheorem{definition}{Definition}[section]
\newtheorem{lemma}[theorem]{Lemma}

\newtheorem{remark}[theorem]{Remark}

\title{A New Approach to Inverse Local Times}
\author{ A.R. Baigorri\thanks{baig@unb.br}\\ Mathematics Department \\ UnB }    
       
\date{July 1st, 2016}

\begin{document}
\maketitle

\begin{abstract}
In 1981 F. Knight published an article with a partial solution to a problem proposed by Ito-McKean \cite{Ito}, [p.217]. 
In this paper Knight \cite{Knight} 
characterized the Levy measures of gap diffusions also known as quasi-diffusions. The proof is very elegant but it uses quite a lot functional
analysis, more specifically spectral Krein Theory.
We present a new proof of  Knight's Theorem, defined at the beginning of the Introduction as well as the new proof of the same theorem refered
as Theorem \ref{theo7} in Section 2.
\end{abstract}

\vspace{20pt}
\section{Introduction}
In 1981 F. Knight published an article with a partial solution to a problem proposed by Ito-McKean \cite{Ito}. 
In this paper Knight \cite{Knight} 
characterized the Levy measures of gap diffusions also known as quasi-diffusions. The proof is very elegant but it uses functional
analysis, more specifically spectral Krein Theory.

The task of re-proving such result with different techniques was sugested by the late Dr. Martin L. Silverstein to whom I am deeply
grateful for his advice and help. The present work provides a constructive proof of the following Knigth's theorem which 
allows a better probablistic  interpretation of every step of the process.

\begin{flushleft}{\bf Theorem\,:}
The Levy measure $n(dy)$ of a persistent gap diffusion on $[0, \infty)$ reflected at $0$  is absolutely continuous with
density 
                 $$n(y) = \int_0^\infty  e^{-y\,x} \mu( dx )$$
with measure $ \mu( dx ) \geq 0 \,\, \text{on} \,\, (0, \infty) $ such that
                 $$\int_0^\infty \frac{1}{x(1+x)} \mu( dx ) < \infty. \,\, \bullet $$
\end{flushleft}           
        
The techniques used are totally different from those previosly used.

In Section 2 we state but rarely prove the material  needed for the proof of  main result refered as Theorem \ref{theo7}
Exact references will be provided in Appendix 3.
Section 2  contains the new proof of Knight's Theorem A as defined at the beginning of the chapter.

\vspace{20pt}
\section{Main Tools}
Firstly, let us define (see Watanabe \cite{Watanabe}) what is called a  gap-diffusion

\begin{definition}\label{inextensible}
Let $0 < l \leq \infty$. A non negative Borel measure $m(dx)$ on $[0, l]$ is called an inextensible measure on $[0,l]$
if there exist a non-negative Borel measure which is finite on each compact set $m'(dx)$ on $[0, l]$ such that, by extending
$m'(dx)$ on $[0, l]$ so that $m'\{ l \}=0$
$$m(dx)  =
      \begin{cases}
            m'(dx)                             & \text{if\,\,\,} l = \infty  \,\,\, \text{or} \,\,\, m'([0, l]) = \infty  \\
            m'(dx) + \infty.\delta_{l}(dx)     & \text{if\,\,\,} l +  m'([0, l] < \infty.  
      \end{cases}$$
where $\delta_l$ is the unit measure at $x = l$.
\end{definition}

Let $m(dx)$ an inextensible measure on $[0, l]$ amd $E_m \in [0, \infty]$ be the support of $m$. Let uss assume, fora
for simplicity, $0 \in E_m$. 

Let $X - (X_t, P_x)_{x \in \mathbb{R}}$ be a one dimensional standard Brownian motion and $\Phi(t, x), (t \geq 0, x \in \mathbb{R})$
be its local time at $x$ 
\begin{center}
$\boxed{
\int_0^l I_A(X_s) \, ds = 2 \int_A \Phi(t, x) \, m(dx), \,\,\, \forall \,\, A \subset {\cal B}(\mathbb{R}).}$
\end{center}

Let
$$ \Theta(t) = \int_{[0, l]} \Phi(t, x) \, m(dx) $$
and 
$\Theta^{-1}(t)$ be the inverse function  of $\Theta(t)$ (inverse local time) and ${\hat X}(t) = X(\Theta^{-1}(t)).$ 
By the general theory of time change, ${\hat X}(t)$ defines a Markov process on $E_m = E_m \cap [0, l]$ whose life time  is 
identified with the first hitting time for $l$.

\begin{definition}\label{gapdiffusion}
The Markov process ${\hat X}(t) = X(\Theta^{-1}(t))$ above mentioned is called the gap-diffusion process corresponding
to the inextensible measure $m(dx)$.
\end{definition}

Let us introduce the following convenient notation:
\begin{itemize}
 \item \,\, $\tau( t, a, w)$ = local time spent at $a$ up to time $t$ by trajectory $w$ starting at $a$
 \item \,\, $\tau(t, a) = \tau(t, a, w)$
 \item \,\, $\tau(t) = \tau(t,0)$
 \item \,\, $\tau(t, w) = \tau(t, 0, w)$
\end{itemize}

\begin{lemma}\label{lemma1}
Let $[W, { \cal B}, {P}_a \,\, a \in E_1 = \{x_0, x_1\}]$ be a gap-diffusion with two points state and let $\{H_0, H_1\}$
be exponential holding times distributions at point $x_0$ and $x_1$ with rates $\{a_0, a_1\}$.
Thus 
$$ P(H_i > t) = \int_t^\infty h_i(l) \, dl = e^{ -a_i t}, \,\,\,\, i= 0, 1; \,\, a_i > 0. $$

\vspace{15pt}
Define $c(t)$= \{number of jumps from $x_0$ to $x_1$ before time $t$\}. 
Then

$$\boxed{ P_0(\tau^{-1}(t,w)-t \in ds) = \sum_{k=0}^\infty P_0(\tau^{-1}(t,w) - t \in ds \, \vert \, c(t)=k) \quad  \bf{( I)} }$$
with 
\begin{eqnarray*}
P_0(c(t)=k) &=& P_0(H_0^1+H_0^2+....+H_0^{k+1} \geq t) - P_0(H_0^1+H_0^2+....+H_0^k \geq t)\\
            &=& \frac{e^{-a_0 t}(a_0 t)^k}{k !}, \,\,\, k \geq 1 \\
\end{eqnarray*}

and $H_0^i$ are such that $H_0^i=H_0, \,\,\, \forall \, 1 \leq i \leq n.$

Therefore 
$$ P_0(\tau^{-1}(t) -t \in ds) = \sum_{k=0}^\infty \frac{e^{-a_0 t}(a_0 t)^k}{k !}  S^{k*}$$
where 
$S$ is the probability distribuition of $\tau^{-1} -t$ conditioned to $c(t) = 1 \,\,\bullet$
\end{lemma}

\vspace{10pt}
{\bf Proof :} Will be omitted.\,\,$\blacklozenge$

\vspace{20pt}
In this particular case $S$ is absolutely continuous with density 
$$S_0(x) = a_1 e ^{-a_1 x}, \,\, x \geq 0$$
but the argument is general and can 
be applied to a quasi-diffusion with any finite state space. Thus we can conclude tha the inverse local time (at 0) for a gap-diffusion. with finite 
state space is a compound Poisson process with the following probability distribution

$$\tau^{-1}(t) - t = \sum_{k=0}^\infty \frac{e^{-a_0 t}(a_0 t)^k}{k !}  S^{k*}$$
$$ \boxed{S(dx) = P_0(\tau^{-1}(t,w) - t \in dx \, \vert \, c(t)= 1) = S(x) dx \quad {\bf (II)}} $$ 
as {\bf (I)} in Lemma \ref{lemma1}.

\vspace{20pt}
\begin{lemma}\label{lemma2}
Let $\big[W, {\cal B}, P_a, a \in \{0,1,2,...,n\} \big]$ be a gap-diffusion with finite state space. 
Let $n(ds)$ be its Levy measure, then
$$ S(x) dx = n(dx)  \,\,\bullet   $$
\end{lemma}

\vspace{10pt}
{ \bf Proof: }
By Theorem \ref{theo12} we have
$$ E_0 \big[ e^{-z \tau^{-1}(t)} \big] = e^{- t[m_0 z + \int_0^\infty ( 1- e^{-zx} ) n(dx) ] } $$
and by Fubini's Theorem 
$$  E_0 \big[ e^{-z \tau^{-1}(t)} \big] = \sum_{k=0}^\infty \frac{e^{-a_0 t}(a_0 t)^k}{k !} \big( T(z)^k \big ) =
    e^{ - a_0 t } e^{a_0 t T(z)} $$ 
where 
$T(z) = \int_0^\infty e^{-zx } S(x) \, dx $ and $S(x)$ as in {\bf (II)} above.    
Therefore 
$$m_0 z + \int_0^\infty (1 - e^{-zx}) n(dx) = a_0 - a_0 T(z), \,\,\,\, z \geq 0 $$
Differentiating with respect to $z$ in both sides
$$ - a_0 T'(z) = m_0 + \int_0^\infty x e^{-zx} n(dx) $$
$$ -a_0 \frac{d}{dz} \int_0^\infty e^{-zx} S(x) dx = m_0 + \int_0^\infty x e^{-zx} n(dx)$$
$$a_0 \int_0^\infty x e^{-zx} S(x) dx = m_0 \int_0^\infty x e^{-zx} n(dx) $$
Therefore
$$ S(x) dx = n(dx) + m_0 \delta(x) $$
but since $n(dx)$ is a measure on $(0, \infty)$
$$ S(x) dx = n(dx) $$
completing the proof of the lemma.\,\,$\blacklozenge$

\vspace{20pt}
Our next task is to study in more detail the probability distribution $S(x)$ as in formula {\bf (II)} for gap-diffusions with n-points 
state space.  

\vspace{10pt}
{\bf Notation :}  We mean by (i) the expression (ii) 
\begin{itemize}
\item  (i) \,\, $\frac{a_1}{b_1 + x} \ominus \frac{a_2}{b_2 + x}, \,\,\,\ a_i \neq 0, \,\, i = 1,2. $ 
\item  (ii)\,\, $\frac{ a_1 }{ b_1 + x - (\frac{ a_2 }{b_2 + x}) }, \,\,\,\ a_i \neq 0, \,\, i = 1,2. $ 
\end{itemize}

\begin{lemma}\label{lemma3}
Let $\big[W, {\cal B}, P_a, a \in E_2 = \{x_0,x_1,x_2\} \big]$ be a quasi-diffusion with 3 points state space and let $\{H_0, H_1, H_2\}$
be the exponential holding times at the points $\{x_0,x_1,x_2\}$ with rates $ \{ a_0, a_1, a_2\}$ respectively.
Then 
$$L(S)(x) = \frac{a_1}{2(x + a_1)} \ominus \frac{a_1 a_2}{x +  a_2}$$
where $S(x)$ is as in {\bf (II)} \,\,$\bullet$ 
\end{lemma}

\vspace{10pt}
{\bf Proof :}
As in Lemma (\ref{lemma1}) let $c(t)$ be the number of jumps made the particle from $x_0$ to $x_1$ up to time $t$.

Once the particle jumps from point $x_0$ it can do any of the following possible movements before it returns to point $x_0$ again. 
We are only interested in the case $c(t) = 1.$ 

\vspace{20pt}
{\bf Possible Trajectories}  
\begin{itemize}
\item {\bf Case (1)\,\,:} \,\, $< x_1  >$
\item {\bf Case (2)\,\,:} \,\, $<  x_1 . x_2  . x_1 >$
\item {\bf Case (3)\,\,:} \,\, $< x_1  .  x_2  .  x_1  . x_2  . x_1 >$ \\
         .................                    
\end{itemize}

{\sl Case (1)} means that the particle jump to point $x_1$ and then jump to point $x_0$. 

{\sl Case (2)} particle jumps to point $x_1$, then to point $x_2$ back to point $x_1$ and then returns to point $x_0$.

Same for {\sl Case (3)} etc., etc.

\vspace{10pt}
Let $p,q$ be such that $p+q=1, \,\,\, p,q \geq 0$ representing the probability of jumping to the right and 
probability of jumping to the left of any given point of the state space.
We are going to assume works for general $p=q=\frac{1}{2}$ in the following computations but the same argument works for 
general $p$ e $q$.

Thus for $p=q=\frac{1}{2}$
$$ S(x) = \frac{h_1(x)}{2} + \frac{ (h_1)^{2*}* h_2 }{4} + \frac{ (h_1)^{3*} * (h_2)^{2*} }{8}  +  .....$$

or $p \neq q$
$$ S(x) = q h_1(x) + p q \Big( (h_1)^{2*}* h_2 \Big)  + p^2 q \Big(  (h_1)^{3*} * (h_2)^{2*} \Big)  +  .....$$

\vspace{10pt}
\begin{eqnarray*}
L(S)(x)&&= \frac{w_1(x)}{2} \Big( 1 + \frac{w_1(x) w_2(x)}{2} + \Big( \frac{w_1(x) w_2(x)}{2} \Big)^2 + .....\Big) \\ 
       &&= \frac{w_1(x)}{2} \Big( \frac{1}{1-\frac{w_1(x) w_2(x)}{2}} \Big )
\end{eqnarray*}

where $$ h_i = a_i e^{-a_i x}, \quad \quad i=0,1,2 $$ 
and $$ w_i(x) = L(h_i)(x) = \frac{a_i}{x + a_i}, \quad \quad i=1,2$$

Thus, for $p=q=\frac{1}{2}$ 
$$L(S)(x) = \frac{a_1}{2(x+a_1)\big(1- \frac{a_1}{2(x+a_1)} \frac{a_2}{x+a_2} \big)}$$

or $p \neq q$ 
$$L(S)(x) = \frac{q a_1}{(x+a_1)\big(1- \frac{p a_1}{(x+a_1)} \frac{a_2}{x+a_2} \big)}$$

Therefore
\begin{eqnarray*}
 L(S)(x) &&= \,\, \frac{a_1}{ 2(x+a_1) - \frac{a_1 a_2}{(x+a_2)} }  \\
         &&= \,\, \frac{a_1}{2(x-a_1)}\ominus\frac{a_1 a_2}{(x - a_2)} 
\end{eqnarray*}

or
$$ L(S)(x) = \,\, \frac{q\,a_1}{(x-a_1)}\ominus\frac{p\,a_1 a_2}{(x - a_2)} $$
completing the proof of the lemma.\,\,$\blacklozenge$

\begin{lemma}\label{lemma4}
\,\, Let $\big[W, {\cal B}, P_a, a \in E_3 = \{x_0,x_1,x_2,x_3\} \big]$ be a quasi-diffusion with 3 points 
state space and let $\{H_0, H_1, H_2, H_3\}$
be the exponential holding times at the points $\{x_0,x_1,x_2,x_3\}$ with rates $ \{ a_0, a_1, a_2, a_3\}$ respectively.
Then 
$$L(S)(x) = \frac{a_1}{2(x + a_1)}\ominus\frac{a_1 a_2}{2(x +  a_2)}\ominus\frac{a_2 a_3}{(x +  a_3)}$$
where $S(x)$ is as in {\bf (II)} \,\,$\bullet$
\end{lemma}

\vspace{10pt}
{\bf Proof :} Same argument as in Lemma \ref{lemma3} we obtain
$$L(S) = \frac{w_1}{2} \Big( \frac{1}{ 1 - \frac{w_1}{2} g_2 } \Big )$$

where 
$$ g_2(x) = \frac{a_2}{2(x+a_2) - \frac{a_2\,a_3}{(x+a_3)}} $$
and
$$w_1(x) = \frac{a_1}{(x+a_1)}.$$

Therefore 
\begin{eqnarray*}
L(S)(x) &&= \frac{a_1}{2(x+a_1)\Big(1 - \frac{a_1}{2(x+a_1)} 
                        \big[ \frac{a_2}{ 2(x+a_2) - \frac{a_2 a_3}{(x+a_3)} } \big] \Big)}  \\  
&&=\frac{a_1}{2(x + a_1) - \frac{a_1 a_2}{ 2(x + a_2) - \frac{a_2 a_3}{(x + a_3)  } } }  \\ 
&&= \frac{a_1}{2(x + a_1)} \ominus \frac{a_1 a_2}{2(x + a_2)} \ominus \frac{a_2 a_3}{2(x + a_3)}\\
\end{eqnarray*}
establishing the lemma.\,\,$\blacklozenge$

\begin{lemma}\label{lemma5}    
Let $\big[W, {\cal B}, P_a, a \in E_N = \{x_0,x_1,x_2,...,x_N\} \big]$ be a 
quasi-diffusion with $(N+1)$ points state space and let $\{H_0, H_1, H_2,.... , H_N \}$
be the exponential holding times at the points $\{x_0,x_1,x_2,...,x_N\}$ with rates $ \{ a_0, a_1, a_2,... ,a_N \}$ respectively.
Let us assume that  
$$ P(H_i > t) = \int_t^\infty h_i(s) \,ds = e^{- a_i t}, \,\,\,\, a_i > 0, \,\, i = 0, 1, ..., N. $$
Then 
$$L(S)(x) 
= \frac{a_1}{2(x + a_1)} \ominus \frac{a_1 a_2}{2(x + a_2)}\ominus .... \ominus\frac{a_{N-1} a_N}{(x + a_N)} $$
where 
$$ S(x) = P_0 ( \tau^{-1}(t, w) - t \in dx \, \vert \, C(t) = 1 ) $$
with $\tau^{-1}(t, w)$ being the inverse locat time for a given finite state space quasi-diffusion \,\,$\bullet$
\end{lemma}

\vspace{15pt}
{\bf Proof :}
We will use an inductive argument over $N$. For

\begin{itemize}
\item  n = 1 $->$ $E_1 = \{x_0,x_1,\}$  \,\, is proved in Lemma \ref{lemma1} 
\item  n = 3 $->$ $E_3 = \{x_0,x_1,x_2,x_3\}$  \,\, is proved in Lemma \ref{lemma3} 
\item  n = 4 $->$ $E_4 = \{x_0,x_1,x_2,x_3,x_4\}$ \,\, is proved in Lemma \ref{lemma4}
\end{itemize}
Suppose true for $(N-1)$, we shall prove for $N$. We know by previous lemma thanks
$$ L(S)(x) = \frac{w_1}{2} \frac{1}{ 1 - \frac{w_1}{2} g_{N-1} }$$
where 
$$ w_1(x) = \frac{a_1}{x + a_1} $$
and
$$ g_{N-1}(x) = \frac{a_2}{2(x+a_2) - \big(\frac{a_2 a_3}{2(x+a_3)}\, \ominus \, ...\,\ominus\,\frac{a_{N-1} a_N}{(x + a_N}\big)}$$
Thus
\begin{eqnarray*}
L(S)(x) &=& \frac{a_1}{2(x+a_1)
            \Big( 1 - \frac{a_1}{2(x-a_1)}
  {\Big[ \frac{a_2}{ 2(x+a_2)\,\ominus\,.....\,\ominus\,\frac{a_{N-1} a_N}{(x + a_N})} \Big] }\Big)}\\
  &=& \frac{a_1}{2(x+a_1)} \ominus \frac{a_1 a_2}{2(x + a_2)} \ominus \frac{a_{N-1}a_N}{(x+a_N)}\\
\end{eqnarray*}
completing the proof of the theorem.\,\,$\blacklozenge$

\vspace{20pt}
Let $B^{+} = \big[ W; {\cal B}; P_a; a \in \mathbb{R}^{+} = \big]$ be a reflected Brownian motion.
Let $D_k = \big[ W^k; {\cal B}^k; {P_a}^k; b \in E^k = {b_0^k, b_1^k, .... , b_{n_k}^k}, b_0^k = 0 \big], \,\, k = 1,2,...$ be 
a quasi-diffusion with finite state space $E^k$, Levy measure $ n_k(dx) $ and speed measure $m_k(dx)$ and let 
$D = \big[W, {\cal B}, P_a, a \in \mathbb{R}^{+} \big]$ be a quasi-diffusion with Levy measure $m(dx)$.

\begin{theorem}\label{theo6}
If for every $D_k, \,\, k = 1,2,...$ as above we have that the Levy measure $n_k(dx)$ of a persistent gap-diffusion (quasi-diffusion)
with finite state space $E^k = \{ b_0^k, b_1^k, ..., b_{n_k}^k \}$ reflected at $0$ as in Theorem \ref{theo8} has the representation
 
\vspace{5pt} 
$$\boxed{ n_k(dy) = \int_0^\infty e^{-yx} \mu_k (dx) \, dy} \,\,\,\, (*)$$
\vspace{5pt}

with measure $\mu_k(dx) \geq 0$ on $(0, \infty$ such that 
 
\vspace{5pt}
$$ \boxed{\int_0^\infty \frac{\mu_k(dx)}{x(1+x)} < \infty.} \,\,\,\, (**)$$
\vspace{5pt}

Then the Levy measure $n(dy)$ of a persistent quasi-diffusin on $[0, \infty)$ reflected at $0$ has the representation $(*)$
satisfying condition $(**)$ $\bullet$
\end{theorem}

\vspace{10pt}
{\bf Proof :} \,\, Let $\dot{\tau}(t,w)$ be the local time spent at $0$ by $D_k$ and $\tau(t,w)$ be the time spent at $0$ 
by the reflecting Brownian motion.
Let $\{m_k\}_{k=1}^\infty$ be the speed measure with finite support converging  monotonically a.e. to $m(dx)$. By Theorem \ref{theo8} we 
know that changing the time of a Brownian motion we can obtain $D_k, \,\,\, k = 1,2,....$.
By [\sl Note A, (A.1)] we have
$$ \dot{\tau}(t,w)(t, \psi, w) = f_k(\tau^{-1}(t, \psi, w)) = \int_0^\infty( \tau( \tau^{-1}(t, \psi,w),w,y)\, m_k(dy) $$
where

$$ f_k(t) = f_k(t,w) = \int_0^\infty \tau(t, w, y) \, m_k(dy) .$$
Then 
$$ \lim_{n \to \infty} \dot{\tau}_n(t, \psi, w) = \dot{\tau}(t, \psi, w)$$
$$ \lim_{n \to \infty} \dot{\tau}_n(t) = \dot{\tau}(t), \,\,\ \text{ monotonically in distribution}$$
$$ \lim_{n \to \infty} E_0 \{ e^{-z \dot{\tau}^{-1}(t)\}} = E_0 \{ e^{-z \dot{\tau}^{-1}(t)\} }, \,\,\, z \geq 0 \eqno{(6.1)} $$

where above limit is monotone.

By Theorem \ref{theo12} we have that for every $z > 0$ 
$$ \lim_{k \to \infty} \int_0^\infty (1- e^{-zy})\, n_k(dy) = \int_0^\infty (1- e^{-zy})\, n(dy) $$
but, since for every $z > 0$
$$\int_0^\infty (1 - e^{-zy})\, n(dy) = \int_ 0^z \big( \int_0^\infty -y \, e^{-sy} \, n(dy) \Big) ds$$
we have for every positive $z$ 
$$ \lim_{k \to \infty} \int_0^z \big( \int_0^\infty -y\,e^{-sy} \, n_k(dy)\big) ds = 
                           \int_0^z \big( \int_0^\infty -y\,e^{-sy} \, n(dy)\big) ds $$
$$ \int_0^z \big( - \lim_{k \to \infty}  \int_0^\infty y\,e^{-sy} \, n_k(dy)\big) ds =
                           \int_0^z - \big( \int_0^\infty y\,e^{-sy} \, n(dy)\big) ds $$
                           
Therefore, by Theorem \ref{theo14} we know that
$$ \lim_{k \to \infty} n_k(dx) = n(dx), \,\,\, \eqno(6.2)$$
vaguely with the dense set of convergence being equal to ${\mathbb{R}}^+= [0, \infty)$.

\vspace{5pt}
Therefore by Lemma \ref{lemma2} and Lemma \ref{lemma5}
$$n_k(dx) = S_k(x)\,dx$$
where
$$S_k(x) = \int_0^\infty e^{-xs} \, d\Phi_k(s), \,\,\, k=1,2,.....$$

By (6.2) and Lemma \ref{lemma2} it is easy to see that $n(dx) = S(x)\, dx$.

\vspace{5pt}
We claim that 
$$S(x) = \int_0^\infty e^{-xs} \, d\Phi(s)$$

Due to (6.2), the fact that $ \{ S_k(x) \}_{k=1}^\infty $ are completely monotonic and the vague convergence, it is true that
the following follows:
$$ \lim_{n \to \infty} \int_z^\infty \int_0^\infty e^{-xs} \, d\Phi_n(s) \, dx = \int_z^\infty S(x) \, dx, \,\,\, z > 0 $$

$$ \lim_{n \to \infty} \int_0^\infty \big( \int_z^\infty e^{-xs} \, dx \big)  \, d\Phi_n(s) = \int_z^\infty S(x) \, dx, \,\,\, z > 0 $$

$$ \lim_{n \to \infty} \int_0^\infty  \frac{1}{s} \,\, e^{-zs} \, d\Phi_n(s) = \int_z^\infty S(x) \, dx, \,\,\, z > 0 $$

By Theorem \ref{theo14}, we can conclude tha exits $d\Phi$ so that 
\begin{eqnarray*}
\int_z^\infty S(x) \, dx &=& \int_0^\infty e^{-zx} \, d\Phi (x) = \int_0^\infty \frac{e^{-zs}}{x} \, [ x \, d\Phi(x)]= \\
                         &=& \int_0^\infty \big( \int_z^\infty e^{-xy} \, dy \big )\, [x\, d\Phi(x)]=\\
                         &=& \int_z^\infty  \big( \int_0^\infty e^{-xy}\, [x\, d\Phi(x)] \, \big) \, dy, \,\,\, z> 0. \\
\end{eqnarray*}

Therefore
$$ S(x) = \int_0^\infty  e^{-xy} \, d\phi(y), \,\, \text{ for almost all }  x > 0 $$
where
$$ d\phi(y) = y\, d\Phi(y)$$
completing in this way the proof of the theorem.\,\,$\blacklozenge$

\begin{theorem}\label{theo7}  
The Levy measure $n(dy)$ of a persistent gap diffusion on $[0, \infty)$ reflected at $0$  is absolutely continuous with
density 
                 $$n(y) = \int_0^\infty  e^{-y\,x} \mu( dx )$$
with measure $ \mu( dx ) \geq 0 \,\, \text{on} \,\, (0, \infty) $ such that
                 $$\int_0^\infty \frac{1}{x(1+x)} \mu( dx ) < \infty \bullet$$ 
\end{theorem}

\vspace{10pt}
{\bf Proof:} By Theorem \ref{theo6} it suffices to prove the theorem for finite state space quasi-diffusions.
By Lemma \ref{lemma5} and Theorem \ref{theo21} plus the adequate interpretation in terms of the probabilistic model we get
$$ L(S)(z) = \int_0^\infty \frac{d\Phi(t)}{(z + t)}, \,\,\ z > 0$$
where $\Phi(t)$ is a real valued positive step function with discontinuities on $[0, \infty)$ and $S$ as in Lemma \ref{lemma2}.

\vspace{5pt}
By Theorem \ref{theo15}
$$L(S) = L\{ \int_0^\infty e^{-tu} \, d\Phi(u) \}$$
and uniqueness of Laplace transform
$$S(x) = \int_0^\infty e^{-xu} d\Phi(u) $$
$$ n(dx) = \big (\int_0^\infty e^{-xu} d\Phi(u) \big) \, dx.$$

Appying Theorem \ref{theo12} E.3  we know that
\begin{eqnarray*}
\int_0^\infty ( 1 - e^{-s})\,n(ds) < \infty &=&\int_0^\infty ( 1 - e^{-s})\, \big( \int_0^\infty e^{-su}\, d\Phi(u) \big) \, ds < \infty \\
&=&\int_0^\infty \big( \int_0^\infty  e^{-su}(1 - e^{-s})\, ds \big) \, d\Phi(u) < \infty \\
&=&\int_0^\infty ( \frac{1}{u} - \frac{1}{u + 1} ) d\phi(u) < \infty\\
&=&\int_0^\infty \frac{1}{u(1+u)} d\phi(u) < \infty. \\
\end{eqnarray*}
completing the proof of the theorem.$\blacklozenge$


\vspace{20pt}
\section{Appendix}

{ \bf NOTE A: \,\, Brownian Motion}.

Let us assume that $Q = [0, \infty)$, \, $m(dx)$ is a measure on $Q$ that is finite om compact 
subintervals and positive on neighborhoods os $0$
Let $B^+$ denote the standard Brownian motion $(B^+ = \vert B(t) \vert).$
Let
$$\tau^+ = \frac{1}{2} \frac{d}{dx} \int_0^t I_{[0,x)}(B^+(s))\, ds$$
denote the local time of $B^+$ (continuous in (t,x) $P_0$ a.s. by Trotter's theorem).

\begin{theorem}\label{theo8}[\cite{Knight}, p.55],\,\cite{Wall}.
 Let $\tau(t)$ the right continuous inverse of the additive functional 
$$A(t) = \int_{0^-}^\infty  \tau^+ (t,x)m(dx)$$
with $\tau(t) = \infty $ for $t \geq A(\infty)$.
Then the process $X_t = B^+(\tau(t))$ with the usual generated $\sigma$-fields and translation operators $\Theta_{\tau(t)}$
defines a Hunt process on closure (support) of $m(dx)$ for the probabilities $P_x of B^+$. It is called tha gap diffusion 
on $[0, \infty)$ (quassi diffusion) with natural scale and speed measure $m(dx)$.\,$\blacklozenge$
\end{theorem}

%

\begin{remark}\label{A} \cite{Grommer}.
We know that a solution od $ \tau'(X(t),\psi) = t$ is 
$$X(t) = {(\tau')}^{-1} (t, \psi) \,\, \text{for every} \,\, \psi$$
If $X(t) = f(\tau^{-1}(t, \psi))$ it is a solution tools
$$ \tau'(f(X(t),\psi) = \tau'(f^{-1}(f(\tau^{-1}(t, \psi))), \psi) = t.$$
Showing that
$$ {\tau'}^{-1}(t, \psi) = f(\tau^{-1}(t, \psi)), \,\,\,\,\,\,\,\, ({\bf A.1}) $$
where
$$ f(t) = \int_Q \tau(t,\psi) \, m(d\psi).$$
\end{remark}


\begin{theorem}\label{theo12} [\cite{Ito},p.30,Note 1].

If $\{Y(t), \, t \geq 0 \}$ is a real valued process with stationary independent and nonnegative increments, 
$Y(0) = 0$ and right continuous paths defined om some probability space $(\Omega^*, \epsilon, P)$ then
$$E[e^{-\lambda Y(t)} ] = e^{-t g(\lambda)}; \,\,\, \lambda < 0, \quad ({\bf A.2})$$
where
$$g(\lambda) = b \lambda + \int_0^\infty (1 - e^{-\lambda t}) \, n(dt), \quad ({\bf A.3})$$
$$ b \geq 0, \,\,\,  n(dt) \geq 0, \,\,\, \int_0^\infty (1- e^{-t})\,n(dt) < \infty,  \quad ({\bf A.4}) $$
$n(dt)$ is the so called Levy measure of the process \,$\bullet$
\end{theorem}

\vspace{10pt}
{\bf NOTE B: \,\, Laplace Transform} 

\begin{theorem}\label{theo14} [\cite{Feller},Vol.II,p.433].

For n=1,2,.... let $U_n$ be a measure with Laplace transform $w_n$. If $w_n(\lambda) \to w(\lambda)$ for $\lambda > a$,
then $w$ is the Laplace Transform of a measure $U$ and $U_n \to U$.
Conversely, $U_n \to U$ and the sequence $\{w_n(a)\}$ is bounded, then $w_n(\lambda) \to w(\lambda)$ for every $\lambda > a$.
Recall that the sequence $\{U_n\}_{n=1}^\infty$ of  measures is said to converge to $U$ if and only if $U_n(I) \to U(I) < \infty$
for every finite interval of continuity.
\end{theorem}

\begin{remark}\label{theo15} [\cite{Wideer},p.334].

If the integral 
$$ f(s) = \int_0^\infty \frac{d\alpha(t)}{s+t} $$
converges, then 
$$ f(s)= \int_{0^+}^\infty e^{-s t} \phi(t) \,dt, \,\,\, s > 0$$
where
$$ \phi(t) = \int_0^\infty e^{- t u} \, d\alpha(u), \,\,\, t > 0 \,\bullet$$
\end{remark}

\newpage
{\bf NOTE C: \,\, Continued Fractions} 

\begin{remark}\label{theo16} [\cite{Jones},p.17].

A continued fraction is an order pair $<<\{a_n\};\{b_n\},\{f_n\}>$ where $a_1,a_2,...$ and $b_0,b_1,.....$
are complex numbers ath all $a_n \neq 0$ and where $\{f_n\}$ is a sequence in the extended complex plane defined as
follows:
$$ f_n = S_n(0), \,\,\,\, n=0,1,2,...({\bf C.1}) $$
where
$$ S_0(w) = s_0(w),\,\, S_n(w) S_{n-1}(s(n(w)); \,\,\, n=1,2,..,$$
$$ s_0(w) = b_0 + w, \,\, s_n(w) = \frac{a_n}{b_n+w}, \,\, n=1,2,..$$
$f_n$ is called the nth-approximant. If $\{a_n\}$ and $\{b_n\}$ are infinite sequences then it is called non terminating
continued fraction, otherwise is called finite or terminating continued fraction. Corresponding to each continued fraction
there are sequences complex$\{A_n\}$ and $\{B_n\}$ defined by the system of second order linear difference equations.
$$ A_{-1},\, A_0 = b_0,\, B_{-1} = 0, B_0 = 1$$
\begin{eqnarray*}
 A_n &=& b_n A_{n-1} + a_n A_{n-2}, \,\,\,\, n=1,2,..... \,\,\,\, ({\bf C.2})\\
 B_n &=& b_n B_{n-1} + a_n B_{n-2}, \,\,\,\, n=1,2,..... \,\,\,\, 
\end{eqnarray*} 
It can be shown by an induction argument that
eqnarray
$$ f_n = \frac{A_n}{B_n}, \,\,\,\,\,\, n=0,1,2,.....$$
$$f_n = b_0 + \frac{a_1}{b_1} \ominus \frac{a_2}{b_2}\ominus ... \ominus\frac{a_n}{b_n} \,\bullet $$
\end{remark}

\begin{theorem}\label{theo18} \cite{Szego}; [\cite{Wideer},p.41].

The following relation holds for any three consecutive orthogona polynonials 
$$ K_n(x) = ( A_n x + B_n ) K_{n-1})x) - C_n(x), \,\,\,\,\,\, ({\bf C.3})$$
Here $A_n , B_n \,\text{and}\, C_n$ are constants, $A_n > 0$ and $C_n > 0$.
If the highest coefficient of $K_n(x)$ is denoted by $a_n$ we have
$$ A_n = \frac{a_n}{a_{n-1}}, C_n = \frac{A_n}{A_{n-1}} = \frac{a_n a_{n-2}}{(a_{n-1})^2}$$
Now let $\{K_n(x)\}$ be the orthonormal set of polynomias associated with the distribution
$d\psi(t)$ on $[0, \infty)$. The recurrence formula (C.3) then suggests the consideration of the continued fraction
$$ \frac{1}{A_1 x + B_1} \ominus \frac{C_2}{A_2 x + B_2} \ominus \frac{C_3}{A_3 x + B_3} \ominus ....({\bf C.4})$$
with $\frac{K_n(x)}{L_n(x)}$ as the nth approximant \,\,\,\, $n=1,2,3....$ \,$\bullet$
\end{theorem}

\begin{theorem}\label{theo19} [\cite{Jones},p.254].

If$\{L_n\}_{n=1}^n$ \,is a sequence of polynomials satusfying a system of three term recurrence relation of the form 
$$ L_0(x) =1, \,\,\, L_1(x) - l_1 + x , \,\,\, L_n(x) (l_n + x)L_n - k_n L_{n-2}(x)  \,\,\,\, ({\bf C.5})$$
with $k_n > 0$, $l_n$ real for $n = 1,2,...$. 
Thenthere exist a $\psi(t)$ real valued, bounded monotone, non decreasing with infinitely many points of increase
such that $\{ L_n \}_{n=1}^n $ is a system of orthogonal polynomials with respect to $\psi(t)$..\,$\bullet$
\end{theorem}

\vspace{10pt}
Given a J-fraction as in Theorem \ref{theo18} (C.4) with the nth-approximant 
$$f_n = \frac{K_n(x)}{L_n(x)}$$
we know that $\{L_n(x)\}_{n=1}^\infty$ satisfies
relation Theorem \ref{theo19} (C.5).
 
\begin{theorem}\label{theo20} \cite{Szego}; [\cite{Wideer},p.54].

The approximants $\frac{K_n(x)}{L_n(x)}$ of the J-continued fraction are determined by the formula 
$$ K_n(x) = c \int_a^b \frac{L_n(x) - L_n(t)}{x-t}\, d\psi(t), \,\,\,\, n=0,1,...\text{({\sl C.6)}} \,\bullet$$
\end{theorem}

{\bf Proof:}
Notice that (C.6) holds also for $n=0$ and $n=1$.
For $n \geq 2$ we have 
\begin{eqnarray*}
  K_n(x) &=& \int_a^b \frac{ L_n(x)-L_n(t)-(l_n + x)\{L_{n-1}(x)-L_{n-1}(t)\}}{(x-t)} - \\
                 & &\phantom{mmmmmmmmmmmmmmmmmm}- \frac{k_n \{ L_{n-2}(x) - L_{n-2}(t) \})}{(x-t)}\,d\psi(t)\\
         &=& \int_a^b \frac{-L_n(t) + (l_n + x) L_{n-1}(t) - k_n L_{n-2}(t)}{(x-t) } \, d\psi(t)\\
         &=& \int_a^b\frac{-L_n(t) + ((l_n + x) + (l_n + t)-(l_n + t))L_{n-1}(t) - k_n L_{n-2}(t)}{(x-t)}\,d\psi(t)\\
         &=& \int_a^b L_{n-1}(t) \frac{(x-t)}{(x-t)}\,d\psi(t).
\end{eqnarray*}
establishing the theorem.\,$\blacklozenge$

\begin{remark}\label{theo21} [\cite{Szego},p.56]; [\cite{Grommer},p.212-238]; [\cite{Jones},p.336-337]. 

For a real J-fraction
$$ \frac{k_1}{l_1 + z} \ominus \frac{k_2}{l_2 + z} \ominus \frac{k_3}{l_3 + z} \ominus..... \,\,({\bf C.7})$$ 
with $k_n > 0 $ and $l_n$ real for $n = 1,2,3,....$, let $K_n(z)$ and $L_n(z)$ be the nth-approximant numerator and denominator
respectively. Then the zeros $\{ x_k\}_{n=1}^n$ of $L_n(z)$ are real and distinct and $\frac{K_n(z)}{L_n(z)}$ 
has a partial fraction decomposition 
$$\frac{K_n(z)}{L_n(z)} = c \sum_{k=1}^n \frac{\lambda_{n_k}}{(z-x_{x_k}}, \, \text{so that} \,\, \lambda_{n_k} > 0;\,k=0,1,2...
                        \,\,\,({\bf C.8})$$
Thus
$$\frac{K_n(z)}{L_n(z)} = \int_{-\infty}^{\infty} \frac{d\Phi_n(t)}{(z-t)}$$
where $\Phi(t)$ is a step function \,$\bullet$

{\bf Proof:} By Lagrange interpolation arguments and since all the zeros of a 
an orthogonal polynomial are distinct and real we have
$$ K_n(x) = C \sum_{k=1} K_n(x_k) \frac{L_n(x)}{L'_n(x) (x - x_k)}$$
By Theorem (\ref{theo20})
$$ \frac{K_n(x_k}{L'_n(x_k)} = \frac{1}{L'_n(x_k)} \int \frac{L_n(t)}{(t-x_k)} \, d\phi(t) = $$
$$ \int \frac{L_n(t)}{L'_n(x_k)(t - x_k)}  =  \lambda _k, \,\,\,\,\, k=1,2,3,...,n$$

$$K_n(x) = C \sum_{k=1}^n  \lambda_k \frac{L_n(x)}{(x - x_k}$$

$$\frac{K_n(x)}{L_n(x)} = C \sum_{k=1}^n  \frac{\lambda_k}{(x - x_k)} \,\,\,\,\,\, ({\bf C.9})$$

where the $\{\lambda_k\}$ are called Christofell numbers; $ \lambda_k > 0 , \,\,\,\, k=1,2,.... .$.

Now ordering $x_k$ according to size and defining the step function $\Phi_n(t)$ by

$$\Phi_n(t) =
      \begin{cases}
            0                              & \text{if\,\,\,} -\infty < t \leq x_1 \\
            \sum_{k=1}^m \lambda_k         & \text{if\,\,\,} x_m < t \leq x_{m+1}, \,\,\,\, m=1,2,...n, \,\,\,\, ({\text \bf C.10}) \\
            k_1                            & \text{if\,\,\,} x_n < t \leq \infty.
      \end{cases}$$

we obtain 
$$ \frac{K_n(z)}{L_n(z)} = \int \frac{d\phi(t)}{(z - t)}, \,\,\,\,\, \text{for every} \,\,\, n.\, \blacklozenge $$
\end{remark}

\vspace{20pt}
\begin{remark}\label{theo22} \cite{Wideer}. 
A necessary and sufficient condition for a rational fraction to have the form
$$ \frac{f_1}{f_0} = \sum_{p=1}^n \frac{L_p}{(z - x_p)}$$
where the $x_p$ are real and distinct and the $L_p$ positive, is tha  $f_1/f_0$  have a continued J-fraction espansion as 
in Remark \ref{theo16} (C.1) in which $a_i$ are positive and $b_i$ are real.
\end{remark}


\vspace{10pt}
\begin{thebibliography}{99}
\bibitem{Blumenthal_1}Blumenthal, R.M. \& Getoor, R.K ``Markov Processes and Potential Theory",
         {\it Academic Press, New York and London}, 1998.
\bibitem{Blumenthal_2}Blumenthal, R.M. \&  Getoor, R.K ``Local Times for Markov Processes",
         {\it Z. wahrscheinlichkeitstheorie, 3, 50-74,}, 1964.
\bibitem{Dynkin}Dynkin, E.B.,''Markov Processes'', {\it Moscow, 1963. English traslation, Springer, Berli},1965.
\bibitem{Feller}Feller, W., ''An Introduction to Probability Theory and its Applications'',
         {\it  Vol. I, 2nd. Ed., and II, Wiley, New York, 1962},1965. 
\bibitem{Grommer}Grommer, J.,''Gauze transcendente Funktrioen mit lauter reelen Nullstellen'',
         {\it J. Reien Angew. Math, 144, 212-238, 1944}, 1944.  
\bibitem{Ito}Ito, K. \&  McKean, H.P., ''Diffusion Processes and Their Sample Paths'',
         {\it Springer-Verlag, Berlin-Heilderberg-New York}, 1965.     
\bibitem{Jones}Jones, W.B. \& Thron, W.J., ''Continued Fractions Analytic Theory and Applications'',
         {\it Addison-Wesley Pub. Co.}, 1980.     
\bibitem{Knight}Knight, F. B., ''Characterization of the Levy Measure of Invers Local Times of Gap Diffusions``,
         {\it Seminar on  Stochastic Processes, E. Cinlar, K. L; Chung, R. K. Getoor, editors. Progres in Probability and Statistics, 
           vol. Birkhauser, Boston}, 1981.
\bibitem{Szego}Szego, G.,''Orthogonal Polynomials'', {\it Colloqium Publications, vol. 23 Amer. Math. Soc. New York}, 1959         
\bibitem{Wall}Wall, H.S., ''Analytic Theory of Continued Fractions'', {\it  Van Nostrand, New York}, 1948.            
\bibitem{Watanabe}Watanabe, S., ''On Time Inversion of One-Dimensional Diffusion Processes'',
         {\it Z. Wahrscheinlichkeitstheorieverw. Gebiete, 31, 115-124}, 1975. 
\bibitem{Wideer} Wideer, D. W., ''The Laplace Transform'', {\it Princeton University Press, Princeton, New Jersey}.  
\end{thebibliography}
\end{document}